\documentclass[final]{article}
\usepackage{mathrsfs}
\usepackage{amsmath}
\def\no{{\nonumber}}
\def\n{\nabla}

\def\Bc{{\cal B}}

\def\n{\nabla}
\def\p{{\partial}}
\def\no{{\nonumber}}

\def\r{{\rightarrow}}

\def\div{{\mbox{div}}}
\def\min{\mbox{min}}
\def\O{\Omega}

\def\D{\Delta}

\def\Bc{{\mathcal B}}

\def\Tc{{\mathcal T}}

\usepackage{amssymb, latexsym, amsmath, epsfig}
\usepackage{amsfonts}
\usepackage{graphics}
\usepackage{amssymb, latexsym, amsmath, epsfig, cleveref}
\usepackage{psfrag,epsfig,color}

\usepackage{mathptmx}
\usepackage{helvet}
\usepackage{courier}
\usepackage{type1cm}

\usepackage{makeidx}         
\usepackage{graphicx}        
\usepackage{multicol}        
\usepackage[bottom]{footmisc}

\usepackage{graphicx}

\title{{\large \bf An optimal nonconforming finite element method for the Stokes equations }
\thanks{email: jiaaanli@gmail.com, supported in part by NSF of China (No. 11771259 and 11371031).}}

\author{Jian Li $^{1,2}$\\
  {\small 1. Department of Mathematics, School of Arts and Sciences, Shaanxi University of Science and Technology,}\\
 {\small Xian 710021, P. R. China}\\
{\small 2.  Institute of Computational Mathematics and its Applications, Baoji University of Arts and Sciences,}\\
 {\small Baoji, 721007, China}\\
}
\date{}
\begin{document}
\maketitle
\begin{abstract}In this paper, we propose and develop an optimal nonconforming finite element
method for the Stokes equations approximated by the Crouzix-Raviart element for velocity
and the continuous linear element for pressure. Previous result in using the stabilization
method for this finite element pair is improved and then proven to be stable. Then, optimal order
error estimate is obtained and numerical results show the accuracy and robustness of the
method.

\end{abstract}


{\bf AMS(MOS) Subject Classification:} 35Q10, 65N30, 76D05
\vspace{0.5cm}
\section{Introduction}

\vspace{0.5cm} Nowadays, finite element methods have become an
important and powerful tool in many scientific and technological
fields. In particular, stable mixed finite element methods are a
fundamental component in search for efficient numerical methods for
solving the Stokes and Navier-Stokes equations governing
incompressible flows
\cite{Brezzi,Brezzi-Fortin,Girault-Raviart,Glowinski,Temam}. For the
incompressible flows, more researches have been directed toward the
compatibility of the component approximations of velocity and
pressure by satisfying the {\it inf-sup} condition in the past
decades, i.e., the stability condition. Some popular finite element
pairs have been constructed for the incompressible Stokes and
Navier-Stokes flows. However, the search for simpler and more
efficient pairs for velocity and pressure approximations is
still attractive and valuable.

Recently, a class of local stabilized mixed finite element methods
have been developed and analyzed for the Stokes and Navier-Stokes
equations approximated by the lower order finite element pairs.
One of them uses the pressure projection method to stabilize
the lower equal-order finite elements (i.e., $P_1-P_1$ or $Q_1-Q_1$).
In practice, this method can also efficiently stabilize the equal-order
conforming finite element pairs $P_r-P_r,~r=1,2$ for the Stokes equations \cite{li1,li2}
and Darcy equations \cite{Zhangxin}. Also, it can be easily promoted
for solving the problems in elasticity, coupling free fluid and porous media system,
fluid-fluid interaction in different media, etc. In \cite{Li-Chen}, a stabilized finite element
method is established for the Stokes equations approximated by
the Crouzix-Raviart element for velocity
and the continuous linear element for pressure (i.e., $NCP_1-P_1$).

Compared with conforming finite element methods, nonconforming
finite element methods for incompressible flows are more popular due
to their simplicity and small support sets of basis functions.
Furthermore, they seem much easier to fulfill the discrete {\it
inf-sup} condition and can
easily relax the high-order continuity requirement for conforming
finite elements. Therefore, in practice, the nonconforming finite
element methods seem superior to the conforming finite element
methods. Based on the above heuristics and some existing result \cite{Lamichhane} ,
we try to optimize the previous method \cite{Li-Chen} and furthermore
establish the weak coercivity, well posedness and optimal estimates of the corresponding system.

As an example, this paper concentrates on a nonconforming finite element
method for the Stokes equations, which uses the nonconforming and conforming piecewise linear polynomial
approximations for velocity and pressure, respectively. This
method is here defined in such a way that it can be easily
generalized to the corresponding nonlinear problem. The present pair is different from the
Crouzeix-Raviart pair \cite{Crouzeix-Raviart}, the $P_2-P_0$ pair \cite{Girault-Raviart},
the MINI-element $P_1b-P_1$ pair \cite{P1b-P1}, the Taylor-Hood
$P_2-P_1$ pair \cite{P2-P1}, and the conforming $P_1-P_1$ pair and nonconforming $NP_1-P_1$ pair
with stabilization based on local Gauss integrations \cite{li2,Li-Chen}.
Also, a better
approximation for the pressure is obtained with the continuous piecewise linear element.
In this paper, it is shown to be stable and optimal match using the nonconforming and conforming piecewise linear polynomial
approximations for velocity and pressure without any stabilization
treatment. It seems more computationally efficient without a loss of accuracy.

The rest of the paper is organized as follows: In the next section, an abstract functional setting for the stationary Stokes
problem is described, along with some useful statements.
Then, in the third section, the weak formulation, stability and
well-posedness are established. Error estimates of
optimal order for the method are derived in section 4.
Finally, numerical experiments are given to show superiority of the
present method for the Stokes equations.

\section{Preliminaries}

This section focus on the stationary Stokes equations with
homogeneous Dirichlet boundary condition. Let the domain $\O$ be a
bounded, convex and open subset of $R^d, d=2,3$ with Lipschitz continuous
boundary $\p\O$. The Stokes equations are presented as follows
\begin{eqnarray}
-\nu\Delta u+\nabla p&=&f \qquad \mbox{in}~~ \Omega,\label{stokes-model-1}\\
\div\ u&=&0 \qquad \mbox{in} ~~\Omega, \label{stokes-model-2}\\
 u&=&0\qquad \mbox{on} ~~\partial\Omega,\label{stokes-model-3}
\end{eqnarray}
where $u=(u_1,u_2,u_d)$ represents the velocity vector, $p$ the
pressure, $f$ the prescribed body force, and $\nu>0$ the viscosity.

For convenience, set
$$
X=\big[H^1_0(\Omega)\big]^d,\quad Y=\big[L^2(\Omega)\big]^d,\quad
M=\left\{q\in L^2(\Omega) : \int_{\Omega}q\ dx=0\right\},
$$$$D(A)=[H^2(\O)]^d\cap X.$$
The spaces $[L^2(\Omega)]^m$, $m=1, 2$, or $4$, are endowed with the
$L^2$-scalar product $(\cdot,\cdot)$ and the $L^2$-norm
$\|\cdot\|_{0}$, as appropriate. The space $X$ is equipped with the
usual scalar product $(\nabla u,\nabla v)$ and the norm
$|\cdot|_{1}$. Note that the norm equivalence between $\|\cdot\|_1$
and $\|\n u\|_0$ on $H^1_0(\O)$, we use the same notation for them.
In fact, standard definitions are used for the Sobolev spaces
$W^{m,r}(\O)$, with the norm $\|\cdot\|_{m,r}$ and the seminorm
$|\cdot|_{m,r}$, $m, r\geq 0$. We will write $H^m(\O)$ for
$W^{m,2}(\O)$ and $\|\cdot\|_{m}$ for $\|\cdot\|_{m,2}$.

Then, the weak formulation of \eqref{stokes-model-1}-\eqref{stokes-model-3} is to seek $(u,p)\in
 X\times M$ such that
\begin{eqnarray}
&B((u,p);(v,q))=(f,v)\qquad\forall (v,q)\in  X\times M,\label{stokes-weakform}
\end{eqnarray}
where
$$B((u,p);(v,q))=a(u,v)-d(v,p)-d(u,q)$$ with
$$
a(u,v)=(-\Delta u,v)=(\n u,\n v),
$$
and
$$ d(v,p)=(\div\ v, p).
$$
Clearly, the bilinear forms $a(\cdot,\cdot)$ and $d(\cdot,\cdot)$
are continuous on $X\times X$ and $X\times M$, respectively. Also,
the bilinear form $d(\cdot,\cdot)$ satisfies the {\it inf-sup}
condition \cite{Girault-Raviart,Temam}
\begin{eqnarray}
\sup_{0\ne v\in X}\frac{|d(v,q)|}{\|v\|_{1}} \geq
\beta\|q\|_{0},\label{inf-sup}
\end{eqnarray}
where $\beta$ is a positive constant depending only on $\O$.

 The
well-posedness of the model problem \eqref{stokes-model-1}-\eqref{stokes-model-3}  follows from the
results of the saddle-point problem. Assume that the domain $\O$ is so regular that
ensures a $H^2$-regularity for the solution of \eqref{stokes-weakform}, namely, the
problem \eqref{stokes-model-1}-\eqref{stokes-model-3} has a unique solution $(u,p)\in D(A)\times H^1(\O)$
satisfying the following a priori estimate
\begin{eqnarray}\|u\|_2+\|p\|_1\leq C\|f\|_0,\label{2.6}
\end{eqnarray}where $C$ is a constant depending on $\O$. Subsequently, the
constant $C>0$ (with or without a subscript) will depend only on the
data $(\nu,~\O,~f)$.

\section{The nonconforming finite element}

Let $\Tc_h$ be a regular triangulation of $\O$ into elements
$\{K_j\}$: $\bar\O=\cup \bar K_j$. Denote a boundary segment and an
interior boundary by $\gamma_j=\p\O\cap\p K_j$ and
$\gamma_{jk}=\gamma_{kj}=\p K_j\cap\p K_k$, respectively. The
centers of $\gamma_j$ and $\gamma_{jk}$ are indicated by $\xi_j$ and
$\xi_{jk}$, respectively. The finite element spaces are the
following nonconforming and conforming finite elements for velocity
and pressure:
$$
\begin{array}{l@{}l}
&\mathbb{NCP}_1=\{v\in Y : v|_K\in\big[P_1(K)\big]^d,~v(\xi_{kj})
=v(\xi_{jk}),~v(\xi_j)=0\quad\forall j,~k,~K\in \Tc_h\},\\[3pt]
&\mathbb{P}_1=\{q\in H^1(\O)\cap M : q|_K\in P_1(K)\quad\forall K\in \Tc_h\}.
\end{array}
$$
We will also use the piecewise constant spaces
$$
\begin{array}{l@{}l}
&\mathbb{P}_0=\{q\in M : q|_K\in P_0(K)\quad\forall K\in \Tc_h\}.
\end{array}
$$
Note that the nonconforming finite element space $\mathbb{NCP}_1$ is not a
subspace of $X$ any more.

Define the energy norm
$$
\|v\|_{1,h}=\left(\sum_{j}|v|_{1,K_j}^2\right)^{1/2}, \quad v\in
\mathbb{NCP}_1.
$$
The two finite element spaces $\mathbb{NCP}_1$ and $\mathbb{P}_1$ satisfy the
approximation property: For $(v,q)\in [H^2(\O)]^d\times H^1(\O)$ there
are two approximations $v_I\in \mathbb{NCP}_1$ and $q_I\in \mathbb{P}_1$ such that
\begin{eqnarray}
\|v-v_I\|_{0}+h(\|v-v_I\|_{1,h}+\|q- q_I\|_{0,})\leq
Ch^2(\|v\|_2+\|q\|_1).\label{interpolation-property}
\end{eqnarray}

Note that the following compatibility conditions hold for all $j$
and $k$:
\begin{eqnarray}\int_{\gamma_{jk}}[v]ds=0\qquad\forall v\in
\mathbb{NCP}_1\label{2.8}
\end{eqnarray}
and
\begin{eqnarray}\int_{\gamma_{j}}vds=0\qquad\forall v\in  \mathbb{NCP}_1,\label{2.9}
\end{eqnarray}
where $[v]=v|_{\gamma_{jk}}-v|_{\gamma_{kj}}$ denotes the jump of
the function $v$ across $\gamma_{jk}$. These conditions also hold
for the rotated $Q_1$ space with the mean integral values as the
degrees of freedom.

\subsection{The weak formulation}

Set $(\cdot,\cdot)_j=(\cdot,\cdot)_{K_j}$,
$\left<\cdot,\cdot\right>_j=(\cdot,\cdot)_{\p K_j}$, and
$|\cdot|_{m,j}=|\cdot|_{m,K_j}$. Then the discrete bilinear forms
are given as follows:
$$
a_h(u,v)=\sum_{j}(\n u,\n v)_j,\quad d_h(v,q)=\sum_{j}(\div\
v,q)_j,$$ $$u|_j,~v|_j\in \left[H^1(K_j)\right]^d,~q\in
L^2(\Omega).\eqno{(2.10)}
$$

Below we mention a few pairs of mixed finite element spaces for the
Stokes equations. Earlier, the lowest-order Crouzeix-Raviart element
using a nonconforming piecewise linear velocity and a piecewise
constant pressure was constructed \cite{Crouzeix-Raviart}, which
was extended to a nonconforming piecewise bilinear velocity and a
piecewise constant pressure in \cite{Cai-Douglas-Ye}.
Recently, we proposed the pressure projection method of noncomforming finite element method
for the Stokes equations approximated by the Crouzeix-Raviart element and
piecewise linear element for velocity and pressure, respectively.

 In general,
the number of degrees of freedom for velocity should be larger than
that for pressure. The $\mathbb{NCP}_1-\mathbb{P}_0$ pair is shown to be stable. Is optimal mixed
finite element space for the possible choice $\mathbb{NCP}_1-\mathbb{P}_1$ stable?
It is still under develop.

As noted, we can not recognize the stable of the finite element pair until
\cite{Lamichhane} have given us a hint related to the choice $\mathbb{NCP}_1-\mathbb{P}_1$. Obviously,
the bilinear form $a_h(\cdot,\cdot)$ and $d_h(\cdot,\cdot)$ are continuous
and coercive with respect to broken-norm:
\begin{eqnarray}|a_h(u_h,v_h)|&\leq& C\|u_h\|_{1,h}\|v_h\|_{1,h},\label{continuouty}\\
a_h(v_h,v_h)&\geq& C\|v_h\|_{1,h}^2.\label{coercivity}
\end{eqnarray}
For completeness, we will prove the inf-sup property which is important for the
incompressible flow.

Finally, the discrete weak formulation of the Stokes equations
(2.1)--(2.3) is to find $(u_h,p_h)\in \mathbb{NCP}_1\times\mathbb{P}_1$ such that
\begin{eqnarray}\Bc_h((u_h,p_h);(v_h,q_h))=(f,v_h)\qquad
\forall (v_h,q_h)\in \mathbb{NCP}_1\times\mathbb{P}_1,\label{weak-formulation}
\end{eqnarray}
where
$$
\Bc_h((u_h,p_h);(v_h,q_h))=a_h(u_h,v_h)-d_h(v_h,p_h)-d_h(u_h,q_h)
$$
is the bilinear form on $(\mathbb{NCP}_1,\mathbb{P}_1)\times (\mathbb{NCP}_1,\mathbb{P}_1)$. In the next
two sections, we will study \eqref{weak-formulation} in terms of stability, and existence and uniqueness.

\subsection{Stability}

In this section, we will study the stability of the nonconforming finite element method \eqref{weak-formulation}
for the Stokes equations. The main result of this subsection is the existence and uniqueness of the nonconforming finite element solution.
The proof of this theorem is based on the inf-sup property that is
proven in the following lemma by adapting to a classical argument \cite{Lamichhane}. For completeness, we provide a detailed proof as follows.

{\bf Theorem 3.1.} There exists a strictly positive constant $\beta>0$ independent of $h$ such that for
every $q_h\in \mathbb{P}_1$, there exists a vector $v_h\in \mathbb{NCP}_1$ such that
\begin{eqnarray}
\sup_{v_h\in X_h}\frac{d_h(v_h,q_h)}{\|v_h\|_{1,h}}&\geq& \beta\|q_h\|_0.\label{inf-sup}
\end{eqnarray}

{\bf Proof.} First, we set a auxiliary space $R_h=\{q_h\in Q_h: q_h=\sum\limits_{i=1}^N q_i\chi_i\}$
where $\chi_i$ is a characteristic function of the support $S_i$ of the standard linear basis function
$\phi_i\in \mathbb{P}_1$ associated with the vertex $x_i, i=1,2,\cdots N$ and $Q_h=\{q_h\in M: q_h|_K\in P_0(K),\int_{\O}q_hdx=0\}$. The relevant finite element space
is also defined as follows
$$M^1_h=\{q_h\in M: q_h|_K\in P_1(K), K\in\Tc_h\}.$$
Furthermore, we define an interpolation operator $I_h: M_h^1\r R_h$ by
$$I_hq_h=\sum_{i=1}^Nq_i\chi_i$$
with
$q_h=\sum\limits_{i=1}^Nq_i\phi_i\in \mathbb{P}_1.$
Then, setting $\div_h$ the divergence operator on each element, we observe that
\begin{eqnarray}
\int_{\O}\div v_hI_hq_hdx&=&\int_{\O}\div v_h\sum_{i=1}^Nq_i\chi_idx\no\\
&=&\sum_{i=1}^Nq_i\int_{S_i}\div v_hdx\no\\
&=&\sum_{i=1}^Nq_i\sum_{K\subset S_i}|K|\div_hv_h.\label{theorem7-1-1}
\end{eqnarray}
On the other hand, letting $P_i$ be the barycentre of the element $K$,
\begin{eqnarray}
d_h(v_h,q_h)&=&\int_{\O}\div v_h\left(\sum_{i=1}^Nq_i\phi_i\right)dx\no\\
&=&\sum_{i=1}^Nq_i\div_hv_h\sum_{K\subset S_i}\int_{K}\phi_idx\no\\
&=&\frac1{d+1}\sum_{i=1}^Nq_i\div_hv_h\sum_{K\subset S_i}|K|,\label{theorem7-1-2}
\end{eqnarray}which together with \eqref{theorem7-1-1} yields
\begin{eqnarray}
d_h(v_h,q_h)=\frac{1}{d+1}\int_{\O}\div_hv_hI_hq_hdx=\frac1{d+1}d_h(v_h,I_hq_h).\label{theorem7-1-3}
\end{eqnarray}
Then, noting that the lower order Crouzeix-Raviart element is stable, namely,
\begin{eqnarray}
\sup_{v_h\in \mathbb{NCP}_1}\frac{d_h(v_h,I_hq_h)}{\|v_h\|_{1,h}}\geq\beta_1\|q_h\|_0,\label{theorem7-1-4}
\end{eqnarray}where $\beta_1>0$ only depends on $\O$, we can obtain that
\begin{eqnarray}
\sup_{v_h\in \mathbb{NCP}_1}\frac{d_h(v_h,q_h)}{\|v_h\|_{1,h}}&=&\sup_{v_h\in \mathbb{NCP}_1}\frac{1}{d+1}\frac{d_h(v_h,I_hq_h)}{\|v_h\|_{1,h}}\no\\
&\geq&\frac{\beta_1}{d+1}\|q_h\|_0\no\\
&=&\beta\|q_h\|_0,\label{theorem7-1-5}
\end{eqnarray}where $\beta=\frac{\beta_1}{d+1}$. $\#$

Furthermore, we can obtain the following result.

{\bf Theorem 3.2.} The bilinear form $\Bc_h((\cdot,\cdot);(\cdot,\cdot))$ satisfies the continuous property
\begin{eqnarray}
|\Bc_h((u_h,p_h);(v_h,q_h))|\leq C(\|u_h\|_{1,h}+\|p_h\|_{0})
(\|v_h\|_{1,h}+\|q_h\|_{0}),\quad (u_h,p_h),~(v_h,q_h)\in
\mathbb{NCP}_1\times \mathbb{P}_1\label{weak-formulation}
\end{eqnarray}
and the coercive property
\begin{eqnarray}
\sup_{0\ne(v_h,q_h)\in\mathbb{NCP}_1\times \mathbb{P}_1}
\frac{|\Bc_h((u_h,p_h);(v_h,q_h))|}{\|v_h\|_{1,h}+\|q_h\|_{0}}
\geq\beta^{*}(\|u_h\|_{1,h}+\|p_h\|_{0}),\quad (u_h,p_h)\in
\mathbb{NCP}_1\times \mathbb{P}_1,\label{weak-formulation-1}
\end{eqnarray}where $\beta>0$ only depends on $\O$.

{\bf Proof.} Using the continuous property of the bilinear forms
$a_h(\cdot,\cdot)$ and $d_h(\cdot,\cdot)$, we can easily obtain the
continuous property of $\Bc_h((\cdot,\cdot);(\cdot,\cdot))$.

As for the weakly coercivity of $\Bc_h((\cdot,\cdot);(\cdot,\cdot))$, there exists a positive constant $C_0$ and $w\in X$ for all $p_h\in
\mathbb{P}_1\subset M$,
such that
\begin{eqnarray} (\div w,p_h)_j=\|p_h\|_{0,j}^2,~~\|w\|_{1,j}\leq
C_0\|p_h\|_{0,j}.\label{weak-formulation-2}
\end{eqnarray}
Setting the finite element approximation $w_h\in X_h$ of $w$, we
have
\begin{eqnarray}\|w_h\|_{1,h}\leq C_1\|p_h\|_0.\label{weak-formulation-3}
\end{eqnarray}
First, taking $(v_h,q_h)=(u_h-\alpha w_h,-p_h)$ for some positive constant $\frac1{\nu C_0^2}>\alpha$ yet to be determined in the bilinear term $\Bc_h((\cdot,\cdot);(\cdot,\cdot))$ to obtain
\begin{eqnarray}
\Bc_h((u_h,p_h);(u_h-\alpha w_h,q_h))&=&a_h(u_h,u_h-\alpha w_h)-d(u_h-\alpha w_h,p_h)+d(u_h,p_h)\no\\
&=&a(u_h,u_h)+\alpha d_h(w_h,p_h)-\alpha a_h(u_h,w_h)\no\\
&=&\nu\|u_h\|_{1,h}^2+\alpha\|p_h\|_0^2-\alpha a_h(u_h,w_h)\no\\
&\geq&\frac{\nu}2\|u_h\|_{1,h}^2+\alpha(1-\frac{\nu C_0^2\alpha}{2})\|p_h\|_0^2\no\\
&\geq&\frac{\nu}2\|u_h\|_{1,h}^2+\frac{\alpha}2\|p_h\|_0^2\label{weak-formulation-4}
\end{eqnarray}
since
\begin{eqnarray}
\alpha a_h(u_h,w_h)&\leq&\nu\alpha\|u_h\|_{1,h}\|w_h\|_{1,h}\leq\nu\alpha C_0\|u_h\|_{1,h}\|p_h\|_0\no\\
&\leq&\frac{\nu}2\|u_h\|_{1,h}^2+\frac{\nu C_0^2\alpha^2}2\|p_h\|_{0}^2.\label{weak-formulation-5}
\end{eqnarray}
Using a triangle inequality to obtain
\begin{eqnarray}
\|u_h-\alpha w_h\|_{1,h}+\|p_h\|_0&\leq&C_3(\|u_h\|_{1,h}+\|p_h\|_0).\label{bound}
\end{eqnarray}
Setting $\beta^{*}=C_2C_3^{-1}$ and choosing $C_2=\min\{\frac{\nu}2,\frac{\alpha}2\}$, we have
\begin{eqnarray}
\Bc_h((u_h,p_h);(u_h-\alpha w_h,q_h))
&\geq&C_2(\|u_h\|_{1,h}^2+\|p_h\|_0^2),\label{}
\end{eqnarray}which together with \eqref{bound} yields the following
\begin{eqnarray}\sup\limits_{(v_h,q_h)\in\mathbb{NCP}_1\times \mathbb{P}_1}\frac{\Bc_h((u_h,p_h);(v_h,q_h))}{\|v_h\|_{1,h}+\|q_h\|_0}&\geq& \frac{\Bc_h((u_h,p_h);(u_h-\alpha w_h,p_h))}{\|u_h-\alpha w_h\|_{1,h}+\|p_h\|_0}\no\\
&\geq &C_2C_3^{-1}(\|u_h\|_{1,h}+\|p_h\|_0)\no\\
&=&\beta^{*}(\|u_h\|_{1,h}+\|p_h\|_0).\label{weak-formulation-6}
\end{eqnarray}

\subsection{The well-posedness}

Based on previous results, we can derive the existence and uniqueness of the nonconforming element solution for the Stokes equations.

{\bf Theorem 5.2.} Under the assumptions of Theorem 3.2, the problem \eqref{weak-formulation} admits a unique solution.

\section{Estimate in energy norm}
In this section, we will derive optimal order error bounds for the Stokes equations. The non-conformity error is controlled as in the following theorem, which is similar to Strange's lemma in nonconforming finite element version for the second order elliptic problem \cite{Ciarlet}. The proof of this lemma requires a bound on the nonconforming error. For completeness, this bound is provided as follows.

{\bf Lemma 4.1.(Stranger's lemma for the Stokes equations)} There exists a constant $C>0$ depending only on the coercivity and the continuity constants such that
\begin{eqnarray}
\sup_{w_h\in \mathbb{NCP}_1}\frac{|a(u,w_h)-d(w_h,p)-(f,w_h)|}{\|w_h\|_{1,h}}&\leq&Ch(\|u\|_2+\|p\|_1).\label{strange's-lemma}
\end{eqnarray}

{\bf Proof.} By the definition of $a(\cdot,\cdot)$ and $d(\cdot,\cdot)$, it follows that

\begin{eqnarray}
a(u,w_h)&=&\sum_j(\n u,\n w_h)_j,~\forall w_h\in \mathbb{NCP}_1\no\\
&=&\sum_j[-(\Delta u,w_h)_j+<\n u,[w_h]\cdot n>_j]\label{strange's-lemma-a}
\end{eqnarray}
and
\begin{eqnarray}
d(w_h,p)&=&\sum_j(\div w_h,p)_j,~\forall w_h\in \mathbb{NCP}_1\no\\
&=&\sum_j[-(\n p,w_h)_j+<p,[w_h]\cdot n>_j].\label{strange's-lemma-d}
\end{eqnarray}Recalling that $\bar{w}_h=\frac1{|s|}\int_sw_hds$ defined above satisfying
\begin{eqnarray}
\int_{\p K}(w_h-\bar{w}_h)ds=0,\label{strange's-lemma-0}\\
\|w_h-\bar{w}_h\|_{0,s}\leq Ch^{1/2}\|w_h\|_{1,K},\label{strange's-lemma-1}
\end{eqnarray}
and noting that a constant and each interior edge appears twice in the
sum of formulation, we can obtain that
\begin{eqnarray}
a(u,w_h)-d(w_h,p)-(f,w_h)&=&\sum_j<\n u+p,[w_h]\cdot n>_j\no\\
&=&\sum_j<\n u+p,([w_h]-[\bar{w}_h])\cdot n>_j.\label{strange's-lemma-equality}
\end{eqnarray}
Recalling the definition of $\mathbb{P}_0^Kw_h=\frac1{|\p K|}\int_{\p K}w_hds$ satisfying
\begin{eqnarray}
\int_{\p K}(w_h-\mathbb{P}_0^Kw_h)ds=0,\label{}\\
\|w_h-\mathbb{P}_0^Kw_h\|_{0,s}\leq Ch^{1/2}\|w_h\|_1,\label{strange's-lemma-2}
\end{eqnarray}we can obtain
\begin{eqnarray}
a(u,w_h)-d(w_h,p)-(f,w_h)&=&\sum_j<\n u+p,[w_h]\cdot n>_j\no\\
&=&\sum_j<(\n u+p)-\mathbb{P}_0^K(\n u+p),([w_h]-\bar{[w_h])}\cdot n>_j\no\\
&\leq&Ch(\|u\|_2+\|p\|_1).\label{strange's-lemma-3}
\end{eqnarray}
Thus, we can achieve the desired result.


{\bf Theorem 4.2.} Under the assumption of Theorems 3.1-3.2, we can obtain that
\begin{eqnarray}
\|u_h-v_h\|_{1,h}+\|p-p_h\|_0&\leq&Ch(\|u\|_2+\|p\|_1).\label{important-result}
\end{eqnarray}

{\bf Proof.} First, multiplying \eqref{stokes-model-1} by $v_h\in X_h$, integrating over $\O$ and applying the Green formula we
have
\begin{eqnarray}
a(u,v_h)-d(v_h,p)-\sum_K<\n u+p,[v_h]\cdot n>=(f,v_h).\label{continous-VEM}
\end{eqnarray}
Using the same approach as for lemma 4.1 and setting $(e_h,\eta_h)=(u_I-u_h,p_I-p_h)$, we find that
\begin{eqnarray}
&&\frac{\Bc_h((e_h,\eta_h);(v_h,q_h))}{\|v_h\|_{1,h}+\|q_h\|_0}\no\\
&=&\frac{\Bc_h((u_I-u,p_I-p);(v_h,q_h))}{\|v_h\|_{1,h}+\|q_h\|_0}+\sum_j<\n u+p,[v_h]\cdot n>_j\no\\
&\leq&C(\|u-u_I\|_{1,h}+\|p-p_I\|_0)+\sum_j<\n u+p,[v_h]\cdot n>_j\no\\
&\leq&Ch(\|u\|_2+\|p\|_1).\label{}
\end{eqnarray}
Then, we have
\begin{eqnarray}
\|e_h\|_{1,h}+\|\eta_h\|_0&\leq&\frac1{\beta_{*}}\frac{\Bc_h((e_h,\eta_h);(v_h,q_h))}{\|v_h\|_{1,h}+\|q_h\|_0}\no\\
&\leq&Ch(\|u\|_2+\|p\|_1).\label{}
\end{eqnarray}
Thus, using a triangle inequality and \eqref{interpolation-property} to obtain \eqref{important-result}. $\#$

\section{Estimate in $L^2$-norm}

The velocity in $L^2$-norm for the nonconforming element method is here analyzed in the same way as it is done for the classical nonconforming methods.
Firstly, we consider the dual problem: Find $(\Phi,\Psi)\in [H^2(\O)\cap X]^d\times L_0^2(\O)$ such that
\begin{eqnarray}
-\Delta \Phi+\n\Psi&=&u-u_h \quad\mbox{ in }~\O,\label{dual-problem-1}\\
\div\ \Phi&=&0 \qquad~~~~\mbox{ in }~\O,\label{dual-problem-2}\\
\Phi|_{\p\O}&=&0 \qquad~~~~\mbox{on}~\p\O.\label{dual-problem-3}
\end{eqnarray}
Because of the convexity of the domain $\O$,
this problem
has a unique solution that satisfies the regularity property
\begin{eqnarray}
\|\Phi\|_{2}+\|\Psi\|_1\leq C\|u-u_h\|_0.\label{regularity}
\end{eqnarray}

{\bf Theorem 5.1.} Under the assumption of Theorems 4.2, we can obtain that
\begin{eqnarray}
\|u-u_h\|_{0}&\leq&Ch^2(\|u\|_2+\|p\|_1).\label{important-result}
\end{eqnarray}

{\bf Proof.} Multiplying \eqref{dual-problem-1} and \eqref{dual-problem-2} by $e=u-u_h$ and $\eta=p-p_h$, respectively,
integrating over $\O$, to obtain that
\begin{eqnarray}
\|e\|_0^2=-\nu\int_{\O}\D\Phi edx+\int_{\O}\n\Psi edx-\int_{\O}\div\Phi\eta dx.\label{L2-estimate}
\end{eqnarray}
Simplified to gives the following
\begin{eqnarray}
\|e\|_0^2&=&a_h(e,\Phi)-d_h(e,\Psi)-d_h(\Phi,\eta)-\sum_j<\frac{\p\Phi}{\p n},e>_j+\sum_j<e\cdot n,\Psi>_j\no\\
&=&a_h(e,\Phi-\Phi_I)-d_h(e,\Psi-\Psi_I)-d_h(\Phi-\Phi_I,\eta)-\sum_j<\frac{\p\Phi}{\p n},e>_j+\sum_j<e\cdot n,\Psi>_j\no\\
&&\quad+a_h(e,\Phi_I)-d_h(e,\Psi_I)-d_h(\Phi_I,\eta).\no
\end{eqnarray}
The difference of \eqref{stokes-model-1} and \eqref{weak-formulation} tested against $v_h=\Phi_I$, implies that
\begin{eqnarray}
a_h(e,\Phi_I)-d_h(e,\Psi_I)-d_h(\Phi_I,\eta)=\sum_j<\n u+p,[\Phi_I]\cdot n>_j.\no
\end{eqnarray}
Thus,
\begin{eqnarray}
\|e\|_0^2&=&a_h(e,\Phi-\Phi_I)-d_h(e,\Psi-\Psi_I)-d_h(\Phi-\Phi_I,\eta)\no\\
&&-\sum_j<\frac{\p\Phi}{\p n},e>_j+\sum_j<e\cdot n,\Psi>_j\no\\
&&+\sum_j<\frac{\p u}{\p n},\Phi_I>_j+\sum_j<u\cdot n,\Phi_I>_j\no\\
&=&E_1+E_2+E_3.\label{L2-estimate-equation}
\end{eqnarray}
Here,
\begin{eqnarray}
|E_1|&\leq&Ch(\|e\|_{1,h}+\|\eta\|_0)(\|\Phi\|_1+\|\Psi\|_0)\no\\
&\leq&Ch^2(\|\Phi\|_2+\|\Psi\|_1)\leq Ch^2\|e\|_0.\label{L2-estimate-equation-1}
\end{eqnarray}
Using the same approach as Lemma 4.1, yields
\begin{eqnarray}
|-\sum_j<\frac{\p\Phi}{\p n},e>_j+\sum_j<e\cdot n,\Psi>_j|&\leq&Ch^2(\|\Phi\|_2+\|\Psi\|_1)\leq Ch^2\|e\|_0,\label{L2-estimate-equation-2}\\
|\sum_j<\frac{\p u}{\p n},\Phi_I>_j+\sum_j<u\cdot n,\Psi_I>_j|&\leq&Ch^2(\|\Phi\|_2+\|\Psi\|_1)\leq Ch^2\|e\|_0.\label{L2-estimate-equation-3}
\end{eqnarray}
Combining \eqref{L2-estimate-equation} with \eqref{L2-estimate-equation-1}-\eqref{L2-estimate-equation-3}, and using \eqref{regularity} and a triangle inequality, yields \eqref{L2-estimate}

\section{Numerical analysis}

This section concentrates on the performance of the nonconforming finite
element method approximated by the Crouzeix-Raviart element and continuous linear element
for the incompressible Stokes equations. We compare the present method with the stable Crouzeix-Raviart element/piecewise
constant element \cite{Ciarlet,Crouzeix-Raviart}, the pressure projection stabilization finite element method approximated by the Crouzeix-Raviart element/continuous
linear element and piecewise linear element/piecewise linear element for the incompressible Stokes equations \cite{Li-Chen,li2}.

In order to illustrate the features of the present method, three test problems are considered to
verify the performance of the present method including a nonphysical example with
a known exact solutions, the driven cavity flow and a flow over a cylinder.

{\bf Problem I(nonphysical example with analytical solution).} In this case, we consider a unit square
with an exact flow solution given by
$$
u(x)=(u_1(x_1,x_2),u_2(x_1,x_2)),~~p(x_1,x_2)=\cos(\pi x_1)\cos(\pi x_2),
$$
$$
u_1(x_1,x_2)=2\pi\sin^2(\pi x_1)\sin(\pi x_2)\cos(\pi x_1),
~~u_2(x_1,x_2)=-2\pi\sin(\pi x_1)\sin(\pi x_2)^2\cos(\pi x_1).
$$ Then, the body force $f(x,t)$ is deduced from the
exact solution and \eqref{stokes-model-1}. We here pay more attention to
convergence rate of four different methods with the same mesh and
the same UMFPACK code. The results in tables 1-4  suggest
that there are no significant differences between three different
nonconforming finite element methods in terms of the relative $H^1$- and $L^2$-norms for velocity.
Obviously, the present method is more efficient than
other methods by comparison. Especially, the $P_1nc-P_1$ and the stabilized $P_1nc-P_1$ schemes have almost achieve the same superconvergence
rate $O(h^2)$ for pressure. However, the latter did not improve on the accuracy of the stailized schemes whilst being significantly more expensive.

\newpage
\mbox{}
\begin{center}  Table~1. The standard Galerkin method for the Crouzeix-Raviart element.
\end{center}
\begin{center}
\begin{tabular}{ccccccc}
\hline
$1/h$&$\frac{\|u-u_h\|_0}{\|u\|_0}$&$\frac{\|u-u_h\|_{1,h}}{\|u\|_{1,h}}$&$\frac{\|p-p_h\|_0}{\|p\|_0}$&$L_2rate$&$H_1rate$&$pL_2rate$\\
  \hline
10&0.151784&0.470627&0.129512&&&\\			
20&0.0426512&0.248179&0.0571999&1.831361482&0.923203046&1.179001248\\
30&0.0194144&0.167253&0.0366254&1.941080342&0.97330847&1.099503145\\
40&0.0110149&0.125928&0.0270155&1.970112883&0.986496373&1.057873405\\
50&0.00707796&0.100926&0.0214341&1.981967072&0.99184018&1.037124587\\
60&0.00492609&0.0841883&0.0177783&1.987917211&0.994570371&1.025685073\\
 \hline
\end{tabular}
\end{center}

\mbox{}
\begin{center}  Table~2. The standard Galerkin method for the $P_1nc-P_1$ pair.
\end{center}
\begin{center}
\begin{tabular}{ccccccc}
\hline
$1/h$&$\frac{\|u-u_h\|_0}{\|u\|_0}$&$\frac{\|u-u_h\|_{1,h}}{\|u\|_{1,h}}$&$\frac{\|p-p_h\|_0}{\|p\|_0}$&$L_2rate$&$H_1rate$&$pL_2rate$\\
  \hline
10&0.152528&0.474482&0.0761601&&&\\		
20&0.0428488&0.250109&0.0224344&1.8317474&0.923796385&1.763322773\\
30&0.0195085&0.168554&0.010413&1.940555072&0.973303601&1.892987718\\
40&0.0110695&0.126911&0.00597363&1.969732367&0.98640182&1.931647106\\
50&0.00711356&0.101716&0.00386744&1.981642554&0.991744713&1.948351144\\
60&0.00495107&0.0848496&0.00270709&1.987692039&0.994420747&1.956535138\\
\hline
\end{tabular}
\end{center}

\mbox{}
\begin{center}Table~3. The stabilized nonconforming finite method for the $P_1nc-P_1$ pair.
\end{center}

\begin{center}
\begin{tabular}{ccccccc}
\hline
$1/h$&$\frac{\|u-u_h\|_0}{\|u\|_0}$&$\frac{\|u-u_h\|_{1,h}}{\|u\|_{1,h}}$&$\frac{\|p-p_h\|_0}{\|p\|_0}$&$L_2rate$&$H_1rate$&$pL_2rate$\\
  \hline
10&0.153018&0.474483&0.0764451&&&\\			
20&0.0429933&0.250109&0.0225213&1.831517616&0.923799426&1.763133924\\
30&0.0195752&0.168554&0.0104535&1.940440287&0.973303601&1.89294877\\
40&0.0111076&0.126911&0.00599673&1.96965315&0.98640182&1.931724578\\
50&0.00713812&0.101716&0.00388229&1.981594896&0.991744713&1.948472776\\
60&0.0049682&0.0848496&0.00271736&1.987652153&0.994420747&1.956786522\\
 \hline
\end{tabular}
\end{center}

\mbox{}
\begin{center}Table~4. The stabilized finite method for the $P_1-P_1$ pair.
\end{center}

\begin{center}
\begin{tabular}{ccccccc}
\hline
$1/h$&$\frac{\|u-u_h\|_0}{\|u\|_0}$&$\frac{\|u-u_h\|_1}{\|u\|_1}$&$\frac{\|p-p_h\|_0}{\|p\|_0}$&$L_2rate$&$H_1rate$&$pL_2rate$\\
  \hline
10&0.0884909&0.269308&0.203893&&&\\			
20&0.0206359&0.130988&0.0961777&2.100372742&1.039822436&1.084037915\\
30&0.00887527&0.0865962&0.0544149&2.080976891&1.020679148&1.404706588\\
40&0.00490909&0.0647394&0.0348421&2.0584535&1.011136883&1.549667603\\
50&0.00311108&0.051703&0.0244392&2.04405932&1.007666018&1.589281981\\
60&0.00214673&0.0430371&0.018301&2.034998961&1.006207251&1.586387822\\
 \hline
\end{tabular}
\end{center}

{\bf Problem II(The driven cavity flow).} The driven cavity is considered for the four different methods.
It is a box full of liquid with its lid moving horizontally at speed one.
The results for both velocity and pressure are
given in Figures 1-2. Numerical result of the present method
shows the same performance as that of other methods.

%

{\bf Problem III(The exterior of a 2d cylinder).} We build a computation mesh the exterior of a 2d cylinder.
 A fluid recirculation zone produced by the hole must be captured correctly.

The geometry for the numerical model of the problem are given in Figure 3. Also, the Diriclet boundary conditions is designed for this model and  $u_1,u_2$ and $p$ denote the velocity components in $x$ and $y$ direction and the pressure. Simulations have been performed with the given
viscosity $\nu=1$. Here, a set of sample results is given in Figure 3. In order to verify the correctness of the method, a comparison of the results
with the standard Taylor-Hood element shows that the present method is creditable.


\end{document}